\documentclass[a4paper, 12pt] {article}
\usepackage[cp1251]{inputenc}
\usepackage[english]{babel}
\usepackage{amsthm}
\usepackage{amsfonts}
\usepackage{amsmath}
\usepackage{mathtext}
\usepackage{cite}
\begin{document}
\begin{center}
{\large\bf On recovering the Sturm--Liouville differential operators on time scales}\\[0.4cm]
{\bf M.A.\,Kuznetsova\footnote{Department of Mathematics, Saratov State University, Astrakhanskaya 83, Saratov,
410012, Russia, email: kuznetsovama@info.sgu.ru, ORCID: https://orcid.org/0000-0003-1083-0799}} \\[0.4cm]
\end{center}

\thispagestyle{empty}

{\bf Abstract.} We study Sturm--Liouville differential operators on the time scales consisting of a finite number of isolated points and segments. In a previous paper it was established that such operators are uniquely determined by their spectral characteristics. In the present paper, an algorithm for their recovery based on the method of spectral mappings is obtained. We also prove that the eigenvalues of two Sturm--Liouville boundary value problems with one common boundary condition alternate.



\section{Introduction}
 We study an inverse spectral problem for the Sturm--Liouville operator on the time scales consisting of a finite number of isolated points and segments.
Inverse spectral problems consist in recovering operators from their spectral characteristics. Such problems have many applications in natural sciences and engineering.  For the classical Sturm--Liouville operator on an interval, inverse problems have been studied fairly completely; the basic results can be found in \cite{[3], [4], [5]}.  

Differential operators on time scales, i.e. closed subsets of the real line, unify classical differential operators and difference operators. 
Namely, they involve the so-called $\Delta$-derivative which may generalize both the classical derivative and the divided difference, depending on time scales structure, see \cite{[1], Hilger90}.
Differential operators on time scales frequently appear in applications, see \cite{[1], finance, fishing}.

Posing and studying inverse spectral problems essentially depend on time scale structure, which causes difficulties for time scales of the general form. The single work devoted to an inverse problem on an arbitrary time scale is \cite{[7]}, where Ambarzumian-type theorem is obtained being an analogue to the simplest result in the  inverse spectral theory for the Sturm--Liouville operator \cite{ambarzumian}. Any further studies always require putting for definiteness some restrictions on the time scale under consideration, see \cite{ozkan2019, yurko2019, asymp, y2-structure}.
For example, in \cite{ozkan2019} a partial inverse problem was studied on the time scale of the form $T=[0, a] \cup T_1,$ where $T_1$ is a time scale  on which the potential was assumed to be given a priori.  
In \cite{y2-structure} recovery of the potential from the spectral characteristics was ivestigated on the time scales consisting of a finite number of segments. In \cite{yurko2019} the case of two segments and a finite number of isolated points between them was studied. The time scales considered in \cite{asymp} generalize those ones in \cite{y2-structure} and \cite{yurko2019}.

In the present paper, as in \cite{asymp}, we consider the Sturm--Liouville operator $\ell$ given by formula \eqref{1} below on quite general class of time scales $T$ consisting of $N < \infty$ segments and $M < \infty$ isolated points:
\begin{equation} \label{time scale}
T = \bigcup_{l=1}^{N+M} [a_l, b_l], \quad  a_{l-1} \le b_{l-1} < a_l \le b_l , \, l = \overline{2, N+M}, \quad a_l < b_l \text{ iff } l \in \{ l_k \}_{k=1}^N,\end{equation}
 where $l_k$ equals to the indice corresponding to the $k$-th segment. 
 In the case when $N=1$ and $M=0$ the operator $\ell$ becomes the classical Sturm--Liouville operator.
   If $T$ consists only of isolated points, i.e. $N=0,$ it is a difference operator.  Inverse spectral problems in the latter case consist in recovering coefficients in the recurrent relations. Such problems were studied in \cite{atkinson, triangle, filomat, gus1, gus2, aktosun, khan} and other works. 
If $N>1,$ then the Sturm--Liouville equation on $T$  is equivalent to $N$ classical Sturm--Liouville equations \eqref{4} subject to matching conditions \eqref{jump conditions}, see the next section. This partially resembles differential equations with discontinuity conditions, matrix equations, and equations on geometrical graphs, see \cite{matrix, matzam, disc, gd, star}.
However, unlike the matching conditions usually imposed along with the mentioned equations, conditions \eqref{jump conditions} possess special dependence on the spectral parameter, which essentially complicates the study.

 In \cite{asymp} we proved that the spectral data uniquely determine the Sturm--Liouville operator  on the time scales of the form \eqref{time scale} having assumed that  $q \in W^1_2[a_{l_k}, b_{l_k}],$ $k=\overline{1,N}.$
 In \cite{y2-structure, yurko2019}  the uniqueness theorems were proved under the weaker  assumption $q(x)\in C(T),$ but the authors were forced to restrict the study to a more particular  structure of $T$ than \eqref{time scale}  and to prespecify $q(x)$ in all isolated points.
 Note that our assumption on $q$ allows to recover it also in the isolated points.

  In this paper, we obtain an algorithm for solving the inverse problem based on the method of spectral mappings.  For this purpose, we need asymptotic formulae proved in \cite{asymp}, which we provide in Section~4  for convenience of the reader. We also prove that the eigenvalues of two boundary value problems with one common boundary condition alternate, see Theorem~1. Our results generalize the well-known ones for the classical continuous and discrete Sturm--Liouville operators.

\section{Sturm--Liouville equation on time scale} \label{Sturm--Liouville equation}

For convenience of the reader, let us first provide some basic notions of time scale theory.
Let  $T$ be so far an arbitrary closed subset of
${\mathbb R},$ which we refer to as a {\it time scale}. Define the so-called {\it jump functions} $\sigma$ and $\sigma_-$ on $T$ in the following
way:
$$
\sigma(x):=\left\{\begin{array}{cl}\inf \{s\in T:\; s>x\}, & x\ne \max T,\\[2mm]
\max T, & x=\max T,
\end{array}\right.
\sigma_-(x):=\left\{\begin{array}{cl}\sup \{s\in T:\; s<x\}, & x\ne \min T,\\[2mm]
\min T, & x=\min T.
\end{array}\right.
$$
A point $x\in T$ is called {\it left-dense}, {\it left-isolated}, {\it right-dense}, and {\it right-isolated}, if $\sigma_-(x)=x,$
$\sigma_-(x)<x,$ $\sigma(x)=x,$ and $\sigma(x)>x,$ respectively. If $\sigma_-(x)<x<\sigma(x),$ then $x$ is called {\it isolated}; if
$\sigma_-(x)=x=\sigma(x),$ then $x$ is called {\it dense}.

Denote $T^0:=T\setminus\{\max T\},$ if $\max T$  exists and is left-isolated, and $T^0:=T,$ otherwise. We also denote by $C(B)$ the class of 
functions continuous on a subset $B \subseteq T.$

A function $f$ on $T$ is called $\Delta$-{\it differentiable} at $t \in T^0,$ if for any $\varepsilon>0$ there exists $\delta>0$ such that
$$
|f(\sigma(t))-f(s)-f^{\Delta}(t)(\sigma(t)-s)|\le\varepsilon |\sigma(t)-s|
$$
for all $s\in (t-\delta, t+\delta)\cap T.$ The value $f^{\Delta}(t)$ is called the $\Delta$-{\it derivative} of the function $f$ at the point~$t.$  

We also introduce derivatives of the higher order $n \ge 2.$ Let the $(n-1)$-th $\Delta$-derivative $f^{\Delta^{n-1}}$ of $f$ be defined on
$T^{0^{n-1}},$ where $a^n = \underbrace{a \ldots a}_{n}$ for any symbol $a.$ If $f^{\Delta^{n-1}},$ in turn, is $\Delta$-differentiable on
$T^{0^n}:= (T^{0^{n-1}})^0,$ then $f^{\Delta^n}:= (f^{\Delta^{n-1}})^{\Delta}$ is called the {\it $n$-th $\Delta$-derivative} of $f$ on
$T^{0^n}.$ For $n\ge1,$ we also denote by $C^{n}(T)$ the class of functions $f$ for which there exists the $n$-th $\Delta$-derivative
$f^{\Delta^n}$ and $f^{\Delta^n} \in C(T^{0^n}).$ 
From now on, $f^{\Delta^\nu}(x_1, \ldots, x_n),$ $\nu \in \mathbb N,$ denotes the $\nu$-th  partial $\Delta$-derivative of the function $f(x_1, \ldots, x_n)$ with respect to the { \it first} argument; we agree that $f^{\Delta^0}(x_1, \ldots, x_n):= f(x_1, \ldots, x_n).$ Further notions of time scale theory can be found in \cite{[1]}.

Consider the following Sturm--Liouville equation on $T:$
\begin{equation}
\ell y:=-y^{\Delta\Delta}(x)+q(x)y(\sigma(x))=\lambda y(\sigma(x)), \quad x\in T^{0^2}, \label{1}
\end{equation}
  where $T$ has form \eqref{time scale}. There $\lambda$ is a spectral parameter, and $q(x)\in C(T^{0^2})$ is a real-valued function. 
Let $N > 0 \text{ or }M \ge 3,$ otherwise Eq.~\eqref{1} degenerates. 

A function $y$ is called a solution of Eq.~\eqref{1}, if $y\in C^2(T)$ and equality \eqref{1} is fulfilled.  For such functions, \eqref{1} is equivalent to the system of equations \begin{equation}
\ell_k y := -y''(x_k)+q(x_k)y(x_k)=\lambda y(x_k),\quad x_k \in (a_{l_k},b_{l_k}),\;k=\overline{1,N},  \label{4}
\end{equation} 
equipped with the {\it jump conditions} 
\begin{equation}
\label{jump conditions} \left.
\begin{array}{cc}
y(a_{l+1}) = \alpha^l_{11}(\lambda) y(b_{l}) + \alpha^l_{12}(\lambda) y^\Delta(b_{l}), & l = \overline{1, N + M -1}, \\[3mm]
y^\Delta(a_{l+1}) = \alpha^l_{21}(\lambda) y(b_{l}) + \alpha^l_{22}(\lambda) y^\Delta(b_{l}), & l =\overline{1,N+M-1-\mu_1},
\end{array}\right\}
\end{equation}
where $$
\begin{array}{cc} \alpha^l_{11}(\lambda) = 1, & \alpha^l_{12}(\lambda) = a_{l+1}-b_l, \\[3mm]
\alpha^l_{21}(\lambda) =  (a_{l+1}-b_l)(q(b_l) - \lambda), & \alpha^l_{22}(\lambda) = 1 + (a_{l+1}-b_l)^2(q(b_l) - \lambda),
 \end{array}$$
 $$ \mu_1 := \delta(a_{N+M}, b_{N+M}), \quad   \delta(k, n) := \left\{ \begin{array}{cc} 1, & k=n,\\ 0, & k\ne n. \end{array}\right.$$
 We arrange the coefficients of the jump conditions into the matrices
\begin{align*} \alpha^l(\lambda) &:= \begin{pmatrix}  \alpha^l_{11}(\lambda) &  \alpha^l_{12}(\lambda) \\
 \alpha^l_{21}(\lambda) &  \alpha^l_{22}(\lambda)\end{pmatrix},\; l \in \overline{1,N+M-1-\mu_1}; \\ 
 \alpha^{N+M-1}(\lambda) &:= \big(\alpha^{N+M-1}_{11}(\lambda),\alpha^{N+M-1}_{12}(\lambda)\big) \; \text{ if }\mu_1 = 1.\end{align*}

Without loss of generality, we assume that $l_{k} < l_{k+1},$ $k=\overline{1, N-1}.$ Denote also  $l_0 := 1,$ $l_{N+1} := N+M,$ $\mu_0 := \delta(a_1, b_1),$  and 
 \begin{equation*} \begin{array}{c}
 \beta^{l_k - s}(\lambda) :=  \alpha^{l_{k}-1}(\lambda) \ldots \alpha^{l_k-s}(\lambda), \quad k=\overline{1, N+\mu_1}, \; s=\overline{1, l_k - l_{k-1}}, \\[4mm]
\beta^{l_N}(\lambda) := (1, 0), \quad  l_N = N+M.
\end{array}
\end{equation*}
By the definition of $\beta^l(\lambda)$ and the jump conditions, we have
\begin{equation} \label{large jumps} \left.
\begin{array}{cc}
\big(y(a_{l_k}), y^\Delta(a_{l_k})\big)^T =  \beta^{l_{k - 1}}(\lambda) \big(y(b_{l_{k-1}}), y^\Delta(b_{l_{k-1}})\big)^T, & k=\overline{2-\mu_0, N}, \\ 
 y(a_{l_{N+1}}) = \beta^{l_{N}}(\lambda) \big(y(b_{l_{N}}), y^\Delta(b_{l_{N}})\big)^T & \text{ if } N+M > l_N, 
 \end{array} \right\}
\end{equation}
where the superscript $T$ is the transposition sign. 

Consider the elements of the matrices $\beta^l(\lambda):$
 \begin{equation*} \beta^l(\lambda) =: \left\{\begin{array}{cc}\begin{pmatrix}  \beta^l_{11}(\lambda) &  \beta^l_{12}(\lambda) \\
 \beta^l_{21}(\lambda) &  \beta^l_{22}(\lambda)\end{pmatrix}, & l=\overline{1, l_N-1}, \\[4mm]
 \big(\beta^{l}_{11}(\lambda),\beta^{l}_{12}(\lambda)\big), &  l=\overline{l_N, l_{N+1}-\mu_1},\end{array}\right. \quad
\end{equation*}
It is easy to see that the coefficients $\beta_{ij}^{l_k}(\lambda)$ are polynomials of the form
\begin{multline} \label{beta polymonials} 
\beta_{ij}^{l_k}(\lambda) = (a_{l_{k-1}+1} - b_{l_{k-1}})^{j} (a_{l_{k}} - b_{l_{k}-1})^{i} \prod_{l=l_{k-1}+1}^{l_{k} - 2} (a_{l+1} - b_l)^2 (-\lambda)^{l_k - l_{k-1} -2+i}  \\
+ O(\lambda^{l_k - l_{k-1} -3+i}), 
\end{multline}
where $k=\overline{2-\mu_0, N+\mu_1},$ $i = \overline{1, 2-\delta(k, N+1)}$ and  $j=1,2.$ The more precise formulae than \eqref{beta polymonials} can be found in \cite[Lemma 1]{asymp}. 

Thus, the study of Eq. \eqref{1} is reduced to the study of system \eqref{4} subject to conditions~\eqref{large jumps} with the coefficients polynomially dependent on $\lambda.$ This allows us to use some auxiliary results from classical Sturm--Liouville theory. 
 \section{Spectral characteristics and inverse spectral problems} \label{Spectral characteristics}
Denote by $L_j(B)$ the boundary value problem for Eq. (\ref{1}) on a closed subset  $B \subseteq T$
with the boundary conditions 
\begin{equation}y^{\Delta^j}(\min B)=y(\max B)=0, \quad j=0,1. \label{bc}\end{equation}
A value $\lambda$ is called an {\it eigenvalue} of the boundary value problem $L_j(B)$ if there exists a non-zero solution $y$ of Eq. \eqref{1} satisfying boundary conditions~\eqref{bc}.

Let $S(x,\lambda)$ and $C(x,\lambda)$ be solutions of Eq. (\ref{1}) on $T$
satisfying the initial conditions
\begin{equation*} \label{initial conditions}
S^{\Delta}(a_1,\lambda)=C(a_1,\lambda)=1,\; S(a_1,\lambda)=C^{\Delta}(a_1,\lambda)=0.
\end{equation*}
For each fixed $x,$ the functions $S(x,\lambda)$ and
$C(x,\lambda)$ are entire in $\lambda.$  We introduce the functions
$$
\Theta_0(\lambda):=S(b_{N+M},\lambda), \ \Theta_1(\lambda) := C(b_{N+M}, \lambda).
$$
For $j=0,1,$ the eigenvalues $\{\lambda_{nj}\}_{n\ge 1}$ of the boundary value problem
$L_j(T)$ coincide with the zeros of the entire function $\Theta_j(\lambda),$ which is called  the {\it characteristic function}
of $L_j(T).$ 

For $k=\overline{1, N},$ denote
\begin{equation*}
d_k := b_{l_k} - a_{l_k}, \quad \delta_k := \delta(l_k, N+M),
\end{equation*}
 $$
 f_{k0}(x) := \left\{ \begin{array}{cc}
\sin d_kx, & \delta_k=1, \\
\cos d_kx, &  \delta_k=0,
\end{array} \right. \quad
 f_{k1}(x) := \left\{ \begin{array}{cc}
\cos d_kx, & \delta_k=1, \\
\sin d_kx, & \delta_k=0.
\end{array} \right.  
$$
In what follows, let us agree that $\sqrt{z}$ is the principal square root of a complex number $z.$ 

In \cite{asymp} we proved the following asymptotic formulae for the characteristic functions:
\begin{equation}\Theta_j(\lambda) = \left\{ \begin{array}{cc}
F_j(\lambda) + O\left(\exp(\gamma_1\tau) \lambda^{N+M-2 +j(1-\mu_0)/2-\mu_1/2} \right), & N> 0,\\[3mm] 
\beta_{1,2-j}^1(\lambda), & N=0,
\end{array}\quad j=0,1,\right. \label{char}\end{equation}
where  
$$\gamma_1 := \sum_{k=1}^N d_k, \quad  \tau := \mathrm{Im}\, \rho \ge 0, \; \rho := \sqrt{\lambda},$$  
$$F_j(\lambda) := (-1)^{j(1-\delta_1)(1-\mu_0)} \rho^{\mu_1+j(1-\mu_0)-1}\prod_{k=1-\mu_0}^{N-1} \beta^{l_k}_{2,2-j\delta(0,k)}(\lambda) \beta^{l_N}_{1, 1+\mu_1}(\lambda) \prod_{k=2}^{N} f_{k0}(\rho) f_{1,(1-\mu_0)j}(\rho).$$

From asymptotics \eqref{char} it follows that $\Theta_j(\lambda)$ are entire functions of order $1/2$ if $N>0;$ otherwise they are polynomials of degree $M-2.$ Then  
Hadamard's factorization theorem gives
$$
\Theta_j(\lambda) = C_j p_j(\lambda), \quad p_j(\lambda) = \lambda^{s_j}\prod_{\lambda_{nj}\ne0} \Big(1 -
\frac{\lambda}{\lambda_{nj}}\Big), \quad j=0,1,
$$
where $C_j$ are non-zero constants, while $s_j$ is the multiplicity of the zero eigenvalue in the spectrum $\{\lambda_{nj}\}_{n\ge
1}.$ By virtue of \eqref{char}, we have
$$
C_j = \lim_{\lambda\to -\infty}\frac{F_j(\lambda)}{p_j(\lambda)}.
$$
Hence, the characteristic functions $\Theta_j(\lambda)$ are uniquely determined by their zeros $\{\lambda_{nj}\}_{n\ge 1}.$ 

 We also introduce the {\it Weyl solution}  $\Phi(x,\lambda)$ as the solution of Eq. (\ref{1}) under the boundary conditions
\begin{equation}
\Phi^\Delta(a_1,\lambda)=1,\quad \Phi(b_{N+M},\lambda)=0.                               \label{8}
\end{equation}

We call $M(\lambda):=\Phi(a_1,\lambda)$ the {\it Weyl function}, which generalizes the classical Weyl function. It is obvious that
\begin{equation} 
\Phi(x,\lambda)=S(x,\lambda)+M(\lambda)C(x,\lambda),   \quad                         
M(\lambda)=-\frac{\Theta_0(\lambda)}{\Theta_1(\lambda)}.                              \label{Weyl}
\end{equation}
In accordance with the second assertion in Proposition~1 below, we put
  $$ \alpha_{n} :=  \mathop{Res}_{\lambda = \lambda_{n1}} M(\lambda) = -\frac{\Theta_0(\lambda_{n1})}{\Theta'_1(\lambda_{n1})}, \quad  n \ge 1.$$
We call $\alpha_n$ a {\it weight number}. The numbers $\{1/\alpha_n\}_{n\ge 1}$ generalize the classical norming constants for the Sturm--Liouville operator (see, e.g. \cite{[5]}). 

 {\it Spectral characteristics} include the Weyl function, two spectra $\{ \lambda_{nj}\}_{n \ge 1},$ and the weight numbers $\{\alpha_n \}_{n \ge 1}.$ The following properties of the eigenvalues and the weight numbers are established in \cite{asymp}. 
\medskip

{\bf Proposition 1.} {\it 1. The spectra $\{\lambda_{n0}\}_{n\ge 1}$ and $\{\lambda_{n1}\}_{n\ge 1}$ have no common elements.

2. All zeros of $\Theta_j(\lambda),$ $j=0,1,$ are real and simple.  

 3.  All weight numbers are positive.
 
 4. For $N>0,$ the spectra $\{\lambda_{n0}\}_{n\ge 1}$ and $\{\lambda_{n1}\}_{n\ge 1}$ are infinite. Otherwise each of them consists of $M-2$ elements. 
 \label{spectral features}
}
\medskip

For $j=0,1,$ assume that eigenvalues $\lambda_{nj}$  are numbered in the increasing order: 
$$\lambda_{nj} < \lambda_{n+1,j}, \quad n=\overline{1, \Theta_{NM}}, \quad \Theta_{NM} := \left\{\begin{array}{cc}
+\infty, & N>0, \\
M-2, & N=0.
\end{array}\right.$$
The following theorem  holds.
\medskip

{\bf Theorem 1.} {\it The eigenvalues alternate in the following way:
$$ \lambda_{n1} < \lambda_{n0} < \lambda_{n+1,1}, \quad n=\overline{1, \Theta_{NM}},$$
where we put $\lambda_{M-1,1} := +\infty$ if $\Theta_{NM} = M-2.$
}

{\it Proof.}
Let us consider another time scale $X = T \bigcup \{b_{N+M} +1\}$ and continue $q$ on $X^{0^2} \setminus T^{0^2}$ arbitrarily. Then the functions $S$ and $C$ can be extended to $X$ to satisfy Eq.~\eqref{1} on $X^{0^2}.$ We also can determine the values $S^\Delta(b_{N+M}, \lambda)$ and $C^\Delta(b_{N+M}, \lambda).$  Let  $\varphi(x, \lambda)$ be a solution of Eq. \eqref{1} on $X^{0^2}$ satisfying the initial conditions $$\varphi(b_{N+M}, \lambda) = 0, \quad  \varphi^\Delta(b_{N+M}, \lambda) = -1.$$
The Wronskian-type determinant $W(y, z):= y(t)z^\Delta(t)-y^\Delta(t)z(t)$ is constant on $X^0$ if $y,z$ obey Eq. \eqref{1} on $X^{0^2}.$ Then $\varphi(a_1, \lambda) = \Theta_0(\lambda)$ and $\varphi^\Delta(a_1, \lambda) = -\Theta_1(\lambda).$

From relation \eqref{1} on $X^{0^2}$ with $y=\varphi$ we get
$$(\lambda-\mu)\int_{a_1}^{b_{N+M}} \varphi(\sigma(t), \lambda) \varphi(\sigma(t), \mu) \Delta t = -\Theta_1(\lambda) \Theta_0(\mu) + \Theta_0(\lambda) \Theta_1(\mu),$$
for details see \cite[Proposition 2]{asymp}.
Dividing on $\lambda-\mu$ both parts and taking the limit at $\mu \to \lambda,$ we obtain
\begin{equation*}\Theta_1(\lambda)\Theta'_0(\lambda) - \Theta'_1(\lambda)\Theta_0(\lambda) = \int_{a_1}^{b_{N+M}} \varphi^2(\sigma(t), \lambda) \Delta t > 0, \quad \lambda \ne \lambda_{n0},\end{equation*}
since $\varphi(a_1, \lambda) \ne 0.$
Then the function $M^{-1}(\lambda) = -\Theta_1(\lambda)/{\Theta_0(\lambda)}$ is increasing on every interval  $(\lambda_{n-1,0}, \lambda_{n0}),$ $n \ge 1,$ where we put $\lambda_{00} := -\infty.$ 

It is easy to see that 
$$\displaystyle \lim_{\lambda \to \lambda_{n0} \pm 0} M^{-1}(\lambda) = \mp \infty, \quad n \ge 1.$$
  Using formulae \eqref{beta polymonials} and \eqref{char}, we obtain
$$ \lim_{\lambda \to -\infty} M^{-1}(\lambda) = \left\{
\begin{array}{cc}
\displaystyle -\frac{1}{a_2-a_1}, & a_1=b_1, \\[4mm]
-\infty, & a_1 < b_1.
\end{array} \right.$$
Then from the monotonicity of $M^{-1}(\lambda)$ it follows that $\lambda_{n1} \in (\lambda_{n-1,0}, \lambda_{n0}),$ $n \ge 1,$ which finishes the proof. \qedsymbol
\medskip

The properties described in Proposition~1 and Theorem~1 are well known for the classical self-adjoint Sturm--Liouville operator, see \cite{[5]}. In \cite{agarval} alternation of two spectra was proved in the case of a general time scale but under different boundary conditions. 
Further properties of the eigenvalues and the weight numbers, namely asymptotic formulae, can be found in Section~4.

From now on, we assume that $q \in W^1_2[a_{l_k}, b_{l_k}],$ $k=\overline{1,N}.$ In \cite{asymp} we proved that {\it spectral data} of three types uniquely determine the potential:
\begin{enumerate}
\item $M(\lambda);$
\item $\{ \lambda_{nj} \}_{n \ge 1},$ $j=0,1;$
\item $\{ \lambda_{n1} \}_{n \ge 1}$ and $\{ \alpha_n\}_{n \ge 1}.$
\end{enumerate}
Given the spectral data of one type, we can recover them of any other one, see \cite{asymp}.  By this reason, it is sufficient to provide only an algorithm solving the inverse problem with $M(\lambda)$ as the input data:
\medskip

{\bf Inverse problem 1.}  Given $M(\lambda),$ find $q(x)$ on $T^{0^2}.$ \medskip

We should note that together with the Weyl function the stucture of $T$ is known. The recovery of $T$ along with the potential requires a separate investigation.
 
 The main result of the paper is an algorithm for solving Inverse problem~1. We give it in recursive style, which means that on the $m$-th step, $m = \overline{1, N+M-\mu_1-1},$ we reduce the recovery of $q$ on $T_m^{0^2}$ to its recovery on $T_{m+1}^{0^2},$ where
\begin{equation} \label{T_m}
T_m = \bigcup_{l=m}^{N+M} [a_l, b_l], \quad m=\overline{1, N+M-\mu_1}.
\end{equation}
The initial step is to recover the potential on the first segment $[a_1, b_1]$  or in the first isolated point $a_1$ of the time scale. Consider the following auxiliary local inverse problem.    
\medskip

{\bf Inverse problem 2.} Given $M(\lambda),$ find $q$ on $[a_{1}, b_{1}].$
\medskip

In Sections~5 and~6, we obtain two algorithms for solving this local inverse problem. There are different approaches according to whether $a _1 < b_1$ or $a_1 = b_1.$ On the interval the potential is recovered by the method of spectral mappings, see Algorithm~1. In isolated point we use asymptotic relations to find $q(a_1),$ see Algorithm~2. 
Further in Section~6, we obtain Algorithm~3 for solution of Inverse problem~1, which is based on the Algorithms for Inverse Problem~2.

\section{Asymptotic formulae} \label{Asymptotic formulae}
In this section, we provide asymptotic formulae  for the eigenvalues and the weight numbers as well as for the functions $C(x, \lambda)$ and $\Phi(x, \lambda),$ which will be used below for solving Inverse problem~2 in the case $a_1 < b_1.$ 

First, we establish asymptotic formulae for $C(x, \lambda)$ and $\Phi(x, \lambda),$ $x \in (a_{1}, b_{1}),$ if $a_1 < b_1.$ Asymptotics for $\Phi(x, \lambda)$ are obtained via decomposition into Birkhoff solutions and then solving the obtained linear system. 
\medskip

{\bf Lemma 1.} {\it 
Let $a_1 < b_1.$ Then for every fixed $\delta > 0$ the estimates 
\begin{equation}  \label{16S}\left.\begin{array}{c}
\displaystyle C(x+a_{1}, \lambda) = \cos \rho x + O\left(\frac{\exp(\tau x)}{\rho}\right), \\[4mm]
 C'(x+a_{1}, \lambda) = -\rho \sin \rho x + O(\exp(\tau x)),
\end{array} \quad x\in(0,d_1], \quad \rho\in G_\delta, \right\}
\end{equation}
\begin{equation}\label{Phi} \left.
\begin{array}{c}
\displaystyle \Phi(x+a_{1},\lambda) = \frac{f_{10}(\frac{x \rho}{d_1}-\rho)}{\rho f_{11}(\rho)}+O\left(\frac{\exp(-\tau x)}{\rho^2} \right),  \\[4mm]
\displaystyle \Phi'(x+a_{1},\lambda) = (-1)^{\delta_1+1}\frac{f_{11}(\frac{x \rho}{d_1}-\rho)}{f_{11}(\rho)}+O\left(\frac{\exp(-\tau x)}{\rho}\right), 
\end{array} \quad x\in[0,d_1), \quad \rho\in G_\delta
\right\}
\end{equation}
hold, where
$$\quad G_\delta := \left\{ \rho: \Big|\rho - \frac{\pi n}{2d_k}\Big| \ge \delta, \; k=\overline{1, N}, \; n \in \mathbb{Z} \right\}.$$
}

{\it Proof.}
 Since  $C(x, \lambda)$ is the cosine-type solution of the first equation in \eqref{4}, equalities \eqref{16S} are obvious. 
 
 Let us prove \eqref{Phi}.
It is known (see, for example, \cite{[5]}) that for $k=\overline{1, N}$ there exists the fundamental system of solutions 
$\{Y_{1k}(x,\rho), Y_{2k}(x,\rho)\},$ $x \in [a_{l_k},b_{l_k}],$ of the $k$-th equation in \eqref{4} having the asymptotics
\begin{equation}\label{12}
Y_{1k}^{(j)}(x+a_{l_k},\rho)=(i\rho)^j\exp(i\rho x)[1],\quad Y_{2k}^{(j)}(x+a_{l_k},\rho)=(-i\rho)^j\exp(-i\rho x)[1],\quad j=0,1,
\end{equation}
 where $[1] := 1 + O(\rho^{-1})$ uniformly for $x\in[0,d_k]$ as $\rho\to\infty.$ 
Expanding $\Phi(x,\lambda)$ for $x \in [a_{l_k},b_{l_k}]$ with respect to the systems  $\{Y_{1k}(x,\rho), Y_{2k}(x,\rho)\},$ $k=\overline{1, N},$ we get
\begin{equation}\label{14}
\Phi(x,\lambda)=A_{2k-1}(\rho) Y_{1k}(x,\rho)+A_{2k}(\rho) Y_{2k}(x,\rho), \quad x \in[a_{l_k},b_{l_k}],\quad k=\overline{1,N},
\end{equation}
where $A:=\big(A_l(\rho)\big)_{l=1}^{2N}$ is a solution of a certain linear system.

To write this system, denote by $D^1_1$ the following matrix:
$$
D^1_1 := \begin{pmatrix}
r_{12} & s_{12}  & 0 & 0 & 0 & 0 & 0 & \ldots & 0 & 0 & 0 & 0 & 0\\

p_{21} & q_{21} & r_{21} & s_{21} & 0 & 0 & 0 & \ldots & 0 & 0 & 0 &  0  & 0 \\
p_{22} & q_{22} & r_{22} & s_{22} & 0 & 0 & 0 & \ldots & 0 & 0 & 0 &  0  & 0 \\

0 & 0 & p_{31} & q_{31} & r_{31} & s_{31} &  0 & \ldots & 0 & 0 & 0 & 0  & 0 \\
0 & 0 & p_{32} & q_{32} & r_{32} & s_{32} &  0 & \ldots & 0 & 0 & 0 & 0  & 0  \\

\vdots & \vdots &\vdots &\vdots & \vdots & \vdots &  \vdots & \ddots & \vdots & \vdots & \vdots & \vdots & \vdots   \\

0 & 0 &0 &0 & 0 & 0 & 0 & \ldots & 0 & p_{N1} & q_{N1} & r_{N1} & s_{N1} \\
0 & 0 &0 &0 & 0 & 0 & 0 & \ldots & 0 & p_{N2} & q_{N2} & r_{N2} & s_{N2} \\
0 & 0 & 0 & 0 & 0 & 0 & 0 & \ldots & 0 & 0 & 0 &  p_{N+1,1}  & q_{N+1,1}
\end{pmatrix}
$$
with the coefficients given for $\nu=1,2$ by the formulae
\begin{equation*}\begin{array}{c}
p_{m\nu} := \beta^{l_{m-1}}_{\nu1}(\lambda)Y_{1,m-1}(b_{l_{m-1}}, \rho) + \beta^{l_{m-1}}_{\nu2}(\lambda)Y'_{1,m-1}(b_{l_{m-1}}, \rho), \quad m=\overline{2, N+1},\\[3mm]
q_{m\nu} := \beta^{l_{m-1}}_{\nu1}(\lambda)Y_{2,m-1}(b_{l_{m-1}}, \rho) +  \beta^{l_{m-1}}_{\nu2}(\lambda)Y'_{2,m-1}(b_{l_{m-1}}, \rho), \quad m=\overline{2, N+1},\\[3mm]
r_{m\nu} := -Y_{1m}^{(\nu-1)}(a_{l_m}, \rho), \quad s_{m\nu} := -Y_{2m}^{(\nu-1)}(a_{l_m}, \rho), \quad m=\overline{1, N},
\end{array}
\end{equation*}
 The other elements are equal to zero.

Substituting \eqref{14} into \eqref{large jumps} and  (\ref{8}), we obtain the linear system
\begin{equation}D^1_1A=(-1,0,\ldots,0)^T
\label{system}
\end{equation} 
with respect to the vector $A.$ Further, using this system, we estimate the coefficients $A_1(\rho)$ and $A_2(\rho).$

Consider the case $N>1.$  Solving system \eqref{system} by Cramer's formulae, we get
\begin{equation}\label{A1A2}
A_{1}(\rho)=\frac{q_{22} \det D_0^{2} - q_{21} \det D_1^{2}}{\det D^1_1}, \quad A_{2}(\rho)=\frac{p_{21} \det D_1^{2} - p_{22}
\det D_0^{2}}{\det D^1_1},
\end{equation}
where $D_j^{2}$ is the submatrix of $D^1_1$ including columns with the numbers $3, 4, \ldots, 2N-1, 2N$ and rows with the numbers $2+j, 4, \ldots, 2N-1,  2N$ for $j=0,1.$

For $j=0,1,$ denote by ${\cal D}_j^{l_{2}}(\lambda)$ the characteristic function of $L_j(T_{l_2}),$ where $T_m$ is determined in \eqref{T_m}.
One can show that 
\begin{equation}\det {D}_1^{1} =  (-2i\rho)^{N}[1] \Theta_1(\lambda), \quad \det {D}_j^{2} = (-1)^{j+1} (-2i\rho)^{N-1}[1]{\cal D}_j^{l_{2}}, \; j=0,1.\label{det rel}\end{equation}  

From \eqref{char} we obtain
$$\Theta_1(\lambda) = (-1)^{1-\delta_1} \rho^{\mu_1}\prod_{k=1}^{N-1} \beta^{l_k}_{2,2}(\lambda) \beta^{l_N}_{1, 1+\mu_1}(\lambda)
 \prod_{k=2}^{N} f_{k0}(\rho) f_{11}(\rho) [1], \quad \rho \in G_\delta.$$
Since ${\cal D}_j^{l_2}$ is the object  that plays for $T_{l_2}$ the same role as $\Theta_j(\lambda)$ does for $T,$ we have  
$${\cal D}_j^{l_2}= (-1)^{j(1-\delta_2)} \rho^{\mu_1+j-1}\prod_{k=2}^{N-1} \beta^{l_k}_{2,2}(\lambda) \beta^{l_N}_{1, 1+\mu_1}(\lambda)
 \prod_{k=3}^{N} f_{k0}(\rho) f_{2j}(\rho) [1], \quad \rho \in G_\delta, \quad j=0,1;$$ 
 one can also see \cite[Lemma 2]{asymp}.
 
 Then from \eqref{det rel} and the subsequent formulae we get
\begin{equation*} \label{est D}
\frac{\det D_0^{2}}{\det D_1^1} = -\frac{1}{2i\rho^2 \sin \rho d_1 \beta_{22}^{l_1}(\lambda)}[1], \quad \frac{\det D_1^{2}}{\det D_1^1} = O\left(\frac{1}{\sin \rho d_1 \beta^{l_1}_{22}(\lambda) \rho}\right), \quad \rho \in G_\delta.\end{equation*}

Using these estimates along with \eqref{beta polymonials}, \eqref{12}, and \eqref{A1A2},
 we obtain
$$
A_{1}(\rho)=\frac{\exp(-i\rho d_1)}{2 \rho \sin \rho d_1}[1], \quad A_{2}(\rho)=\frac{\exp(i\rho d_1)}{2 \rho \sin \rho d_1}[1], \quad \rho \in G_\delta.
$$
These formulae along with \eqref{14} for $k=1$ yield \eqref{Phi} for $N>1.$ The case $N=1$ can be treated analogously. \qedsymbol
\medskip

In order to write the asymptotic formulae for the eigenvalues and the weight numbers,  we introduce the constants $\delta_k^j := \frac12 \delta(\delta_k, j)$  and
\begin{equation}c_k := \frac12 \int_{a_{l_k}}^{b_{l_k}} q(t) \, dt + \sum_{l=l_k+1}^{\min(l_k+1, N+M)} (a_{l} - b_{l-1})^{-1}, \quad z_k :=\frac{1}{\pi} \Big(c_k + \sum_{l=max(1, l_k-1)}^{l_k-1} (a_{l+1} - b_{l})^{-1}\Big)
\label{z_k}\end{equation}
for $k=\overline{1, N},$ $j=0,1.$

The following two theorems  are Theorems~3 and~4, respectively, from \cite{asymp}. From now on, $\{ \kappa_n\}_{n \ge 1}$ denotes different sequences from $l_2.$
\medskip

{\bf Theorem 2.} {\it 
Fix $j\in\{0,1\}.$ Then the spectrum of $L_j(T)$ consists of $N+1$ parts:
\begin{equation}\{ \lambda_{nj} \}_{n \ge 1} = \Lambda_j \bigcup  \bigcup_{k=1}^N \big\{ (\rho^{k}_{nj})^2 \big\}_{n \ge 1}, \quad j=0,1,
\label{union}
\end{equation}
where $\Lambda_j$ contains $N+M-1+j(1-\mu_0) \mathrm{sign}(N-1+\mu_1)-\mu_1$ elements. Assume that 
\begin{equation} \label{commensurability}
d_k = r x_k, \; x_k \in \mathbb{N}, \quad k = \overline{1, N}, \text{ for some } r \text{ independent of } k.
\end{equation} Then for the sequences $\big\{ \rho^{k}_{nj} \big\}_{n \ge 1}$ the following asymptotic formulae are fulfilled: 
\begin{equation}\rho_{nj}^{k} = \frac{\pi( n-\delta^{j \delta(1, l_k)}_k)}{d_k} + \frac{z_k }{n-\delta^{j \delta(1, l_k)}_k} + \frac{\kappa_n}{n}, \quad k = \overline{1, N}, \ n \in \mathbb{N},
 \label{kappa/n^2}\end{equation}
  where $z_k$ are real constants given in \eqref{z_k}.
 }
\medskip

{\bf Theorem 3.} {\it 
In accordance with \eqref{union}, the sequence $\{ \alpha_n \}_{n \ge 1}$ consists of $N+1$ parts:
\begin{equation*}\label{wn parts}\{ \alpha_n \}_{n \ge 1} = A 
\bigcup \bigcup_{k=1}^N \big\{ \alpha^{k}_n \big\}_{n \ge 1}, \quad
 \alpha^{k}_n := \mathop{Res}_{\lambda = (\rho^{k}_{n1})^2} M(\lambda), \ A := \left\{\mathop{Res}_{\lambda = z} M(\lambda) \colon z \in \Lambda_1 \right\}.\end{equation*}
If \eqref{commensurability} holds, and 
\begin{equation} \label{z restriction}
\frac{z_l}{d_l} \ne \frac{z_\nu}{d_\nu} \text{  if  } l\ne \nu,\; l,\nu=\overline{1,N},
\end{equation} 
then the following asymptotic formulae are fulfilled:
  \begin{equation} \label{wn precise}
  \alpha^{k}_n =\left\{ \begin{array}{cc} \displaystyle\frac{2}{d_{k}} \Big(1  + \frac{\kappa_n}{n}\Big), & k=1, \, a_1\ne b_1, \\[4mm]
 \displaystyle  \frac{\kappa_n}{n}, & \text{ otherwise}.\end{array} \right.
  \end{equation} }
  \medskip
  
We need asymptotic formulae \eqref{kappa/n^2} and \eqref{wn precise} to apply the method of spectral mappings in the following section.
\medskip

{\bf Remark 1.}
Further, we consider implicitly the eigenvalues and the weight numbers of the boundary value problems $L_j(T_{l_k}),$ $k = \overline{1, N}.$ For these boundary value problems, \eqref{commensurability} and the following condition guarantee the analogs of asymptotic formulae \eqref{kappa/n^2} and \eqref{wn precise} with residual summands $\kappa_n/n:$
\begin{equation} \label{condition on z}
\frac{c_k}{\pi d_k} \notin \Big\{ \frac{z_s}{d_s}\Big\}_{s=k+1}^N 
\text{ and } \frac{z_l}{d_l} \ne \frac{z_\nu}{d_\nu} \text{  if  } l \ne \nu,\; l,\nu=\overline{k+1,N}, \text{ for each } k \in \overline{1, N-1}.
\end{equation}


\section{Solution of Inverse problem~2: the case $a_1<b_1$} \label{segment}

In this section, we assume that $a_1 < b_1.$  Let $M(\lambda)$ and $T$ be given. We explain how $q$ can be found on $[a_{1}, b_{1}]$ by the method of spectral mappings. 
For simplicity, we assume that \eqref{commensurability} and \eqref{z restriction} are fulfilled. 
These conditions allow us to use the simplest form of the method of spectral mappings, see \cite[Ch. 1]{[5]}.

Together with the boundary value problem $L_1(T)$ we consider a problem $\tilde L_1(T)$ of the same form but with
another potential $\tilde q.$ If an object $\gamma$ is related to $L_1(T),$ we denote by $\tilde \gamma$ the analogous object for $\tilde L_1(T).$ Put
$$
\xi_n:=\Big|\sqrt{\lambda_{n1}}-\sqrt{\tilde\lambda_{n1}}\Big|+|\alpha_n-\tilde\alpha_n|, \quad n\in {\mathbb N}.$$
It follows from  Theorems~2 and~3 that $\{\xi_n\}_{n\in{\mathbb N}}\in l_2.$

We introduce for $x \in [a_1, b_1]$ the function 
\begin{equation}
D(x,\lambda,\mu):=\displaystyle\frac{\langle C(x,\lambda),
C(x,\mu)\rangle}{\lambda-\mu}=
\displaystyle\int^x_{a_1} C(t,\lambda)C(t,\mu)\, dt,             \label{1.6.2}
\end{equation}
where $\langle y,z\rangle :=yz'-y' z$ and the classical derivatives are taken with respect to the first argument.
The second equality in formula \eqref{1.6.2} follows from the relations $\ell_1 C(t,\lambda)= \lambda C(t,\lambda)$ and $\ell_1 C(t,\mu)= \mu C(t,\mu)$ on $[a_1, b_1]$ (see \eqref{4} for the definition of $\ell_k$).  

Introduce the following designations:
$$
\theta_{n0}:=\lambda_{n1},\;\theta_{n1}:=\tilde{\lambda}_{n1}, \;
\alpha_{n0}:=\alpha_n,\;\alpha_{n1}:=\tilde{\alpha}_n,\;
C_{ni}(x):=C(x,\theta_{ni}),\;\tilde{C}_{ni}(x):=
\tilde{C}(x,\theta_{ni}),
$$
$$
R_{ni,kj}(x) :={\alpha_{kj}}D(x,\theta_{ni},\theta_{kj}),
\; \tilde{R}_{ni,kj}(x):={\alpha_{kj}}\tilde D(x,
\theta_{ni},\theta_{kj}),\; i,j=0,1,\; n,k \in \mathbb{N}.   
$$
Then, in particular, we have $\xi_n = |\sqrt{\theta_{n0}}-\sqrt{\theta_{n1}}|+|\alpha_{n0}-\alpha_{n1}|.$ 

In order to obtain the solution of  Inverse problem~2,  we need the following lemmas. From now on, $C$ denotes sufficiently large positive constants. 
\medskip

{\bf Lemma 2.} {\it The following estimates are valid for
$x\in [a_1, b_1],\; n,k \in \mathbb{N},\; i,j=0,1:$
\begin{equation}
|C_{ni}(x)| \le C ,\quad
|C_{n0}(x)-
C_{n1}(x)|\le C\xi_n,               \label{1.6.6}
\end{equation}
\begin{equation}
\left.\begin{array}{c}
|R_{ni,kj}(x)| \le \displaystyle \frac{C}{|\sqrt{\theta_{n1}}-\sqrt{\theta_{k1}}|+1}, \\[4mm]
|R_{ni,k0}(x)-R_{ni,k1}(x)| \le \displaystyle \frac{C\xi_k}{|\sqrt{\theta_{n1}}-\sqrt{\theta_{k1}}|+1}, \\[4mm]
|R_{n0,kj}(x)-R_{n1,kj}(x)| \le \displaystyle \frac{C \xi_n}{|\sqrt{\theta_{n1}}-\sqrt{\theta_{k1}}|+1}, \\[4mm]
|R_{n0,k0}(x)-R_{n1,k0}(x)-R_{n0,k1}(x)+R_{n1,k1}(x)|\le
\displaystyle \frac{C\xi_n\xi_k}{|\sqrt{\theta_{n1}}-\sqrt{\theta_{k1}}|+1}.
\end{array}\right\}                                            \label{1.6.7}
\end{equation}
The analogous estimates are also valid for
$\tilde{C}_{ni}(x), \; \tilde R_{ni,kj}(x).$ }
\medskip

The proof of this lemma is standard, see \cite[Lemma 1.6.2]{[5]}. In the proof, we also obtained the inequalities
\begin{equation} \label{D est} \left.\begin{array}{c}
|D(x,\lambda,\theta)|\le \displaystyle \frac{C\exp(\tau(x-a_1))}{|\rho - \sqrt \theta|+1}, \quad |Im\,\theta|\le r,\\[4mm]
 |D(x,\lambda,\theta_{k1})-D(x,\lambda,\theta_{k0})| \le \displaystyle
\frac{C\xi_k\exp(\tau(x-a_1))}{|\rho - \sqrt{\theta_{k1}}|+1}.
\end{array}\right\}
\end{equation}
\medskip

{\bf Lemma 3.} {\it The following quantities are finite:}
\begin{equation} \label{ser1}
\sum_{k=1}^\infty \frac{1}{(|\rho - \sqrt{\theta_{k1}}|+1)^2} \le \infty, \quad \lambda \in \mathbb C,
\end{equation}
\begin{equation} \label{ser2}
\sup_{n \in \mathbb N} \, \sum_{k=1}^\infty \frac{1}{(|\sqrt{\theta_{n1}} - \sqrt{\theta_{k1}}|+1)^2} \le \infty.
\end{equation}

{\it Proof.} Let us prove \eqref{ser2}; for \eqref{ser1} the proof is analogous. For simplicity, we assume that all $\theta_{n1} \ge 0,$ which can be achieved by shifting the potential $\tilde q$ in Eq.~\eqref{1}.

Consider the numbers
\begin{equation} \label{R_B}
R_B := \pi^2\Big( \frac{B}{2D}+\frac{1}{4D}\Big)^2, \; B \in \mathbb{N}, 
\quad R_0 := 0, \quad D := x_1 \ldots x_N \, r,
\end{equation}
 where $x_j$ and $r$ are determined in \eqref{commensurability}. By formula \eqref{union} and \eqref{kappa/n^2}, the set $J_B := \{ n \in \mathbb N \colon R_{B-1} \le \theta_{n1} < R_B\}$ contains no more than $N$ elements for a suffiently large $B.$ Then the number of elements in each $J_B$ is bounded by some constant $K.$ Consequently, for any $n, j \in \mathbb N$ we have the inequality
 \begin{equation} \label{KK}
|\sqrt{\theta_{n1}} - \sqrt{\theta_{m1}}| \ge \frac{j}{2D}, \quad 
m \in \mathbb N \colon j+1 \le\frac{|m -n|}{K}< j+2.
\end{equation}
 
 For $n \in \mathbb N,$ we can write
 \begin{multline*}
 \sum_{k=1}^\infty \frac{1}{(|\sqrt{\theta_{n1}} - \sqrt{\theta_{k1}}|+1)^2} = 
 \sum_{k=1}^{n-K} \frac{1}{(|\sqrt{\theta_{n1}} - \sqrt{\theta_{k1}}|+1)^2}\\
  +  \sum_{k=\max(1, n-K+1)}^{n+K-1} \frac{1}{(|\sqrt{\theta_{n1}} - \sqrt{\theta_{k1}}|+1)^2}+
 \sum_{k=n+K}^\infty \frac{1}{(\sqrt{|\theta_{k1}} - \sqrt{\theta_{n1}}|+1)^2}.
 \end{multline*}
 Using inequality \eqref{KK}, we obtain
 $$\sum_{k=1}^\infty \frac{1}{(|\sqrt{\theta_{n1}} - \sqrt{\theta_{k1}}|+1)^2} \le K\sum_{j=1}^{\lceil (n-K)/K\rceil} \frac{1}{(\frac{j}{2D}+1)^2}+(2K-1)+K\sum_{j=1}^{\infty} \frac{1}{(\frac{j}{2D}+1)^2},$$
 where $\lceil \cdot \rceil$ is the ceiling function.
 This concludes the proof. \qedsymbol
\medskip

Let $P(x,\lambda)=[P_{jk}(x,\lambda)]_{j,k=1,2}$ be the matrix determined by the formula
\begin{equation*}
P(x,\lambda) \left[ \begin{array}{ll}
\tilde{C}(x,\lambda) & \tilde{\Phi}(x, \lambda)\\
\tilde{C}'(x,\lambda) & \tilde{\Phi}'(x,\lambda) \end{array}\right]
= \left[ \begin{array}{ll}
C(x, \lambda) & \Phi(x, \lambda)\\
C'(x,\lambda) & \Phi'(x,\lambda) \end{array}\right], \quad x \in [a_{1}, b_{1}].    \label{1.4.15}
\end{equation*}
Since $\langle C, \Phi \rangle = 1,$ the matrix $P(x,\lambda)$ exists and the following formulae hold:
\begin{equation*}
\begin{array}{c}
P_{j1}(x,\lambda)=C^{(j-1)}(x,\lambda)\tilde{\Phi}'(x,\lambda)-
\Phi^{(j-1)}(x,\lambda)\tilde{C}'(x,\lambda),\\[3mm]
P_{j2}(x,\lambda)=\Phi^{(j-1)}(x,\lambda)\tilde{C}(x,\lambda)-
C^{(j-1)}(x,\lambda)\tilde{\Phi}(x,\lambda).
\end{array}                                         \label{1.4.16}
\end{equation*}
Using these formulae for $P_{ij}(x, \lambda)$ along with \eqref{16S} and \eqref{Phi}, we obtain the estimates
\begin{equation}
|P_{ij}(x,\lambda)-\delta(i, j)| \le C|\rho|^{i-j-1},
\quad\rho\in G_{\delta}, \; x \in [a_1, b_1], \; i, j=1,2.               \label{1.4.18}
\end{equation}

Fix $\delta > 0$ and $x \in [a_1, b_1].$   
In the $\lambda$-plane we consider closed contours $\gamma^0_B,$ $B \in {\mathbb N},$
(with counterclockwise circuit) of the form 
 $\gamma_B^0 =\gamma_B^+\cup\gamma_B^-\cup\gamma'\cup(\Gamma_B
\setminus\Gamma_B'),$ where
$$
\gamma_B^{\pm}:=\Big\{\lambda:\; \pm Im\,\lambda=\delta,\;Re\,\lambda
\ge\theta',\;|\lambda|\le R_B\Big\},\quad \theta':=\displaystyle\min_{n\in \mathbb{N},\, i=0,1} \theta_{ni},
$$
$$
\gamma':=\Big\{\lambda:\;\lambda-\theta'=\delta\exp(i\alpha),\;
\alpha\in \Big(\displaystyle\frac{\pi}{2},\displaystyle\frac{3\pi}{2}\Big)\Big\},
$$
$$
\Gamma_B':=\Gamma_B\cap\{\lambda:\;|Im\,\lambda|\le\delta,\;
Re\,\lambda >\theta'\}, \quad \Gamma_B:=\Big\{\lambda:\; |\lambda|= R_B
\Big\},
$$
while $R_B$ are determined in \eqref{R_B}. 
Applying Cauchy's integral formula on the contours $\gamma^0_B$ for the elements of $P(x, \lambda),$ 
analogously to the case of the classical Sturm--Liouville equation \cite[Sect.~1.6]{[5]}
 one can prove the following relations: 
\begin{equation}
\tilde{C}(x,\lambda)=C(x,\lambda)+\displaystyle\sum_{k=1}^{\infty}
\Big(\alpha_{k0}\tilde D(x,\lambda,\theta_{k0})
C_{k0}(x) - \alpha_{k1}\tilde D(x,\lambda,\theta_{k1})
C_{k1}(x)\Big),      \label{1.6.10}
\end{equation}
\begin{equation}
D(x,\lambda,\mu) - \tilde D(x,\lambda,\mu) + \displaystyle\sum_{k=1}^{\infty}\Big(
\alpha_{k0}{\tilde D(x,\lambda,\theta_{k0})}D(x,\theta_{k0},\mu)
-\alpha_{k1}{\tilde D(x,\lambda,\theta_{k1})}D(x,\theta_{k1},\mu)\Big) = 0,          \label{1.6.11}
\end{equation}
where the series converge absolutely and uniformly with respect to
$x\in [a_{1}, b_{1}]$ and $\lambda,\,\mu$ on compact sets. Absolute and uniform convergence follows from estimates \eqref{D est} and \eqref{ser1}.
 
It follows from the definition of $\tilde{R}_{ni,kj}(x),\;
{R}_{ni,kj}(x)$ and formulae (\ref{1.6.10}), (\ref{1.6.11}) that
\begin{equation}
\tilde{C}_{ni}(x)=C_{ni}(x)+\displaystyle\sum_{k=1}^{\infty}
\big(\tilde{R}_{ni,k0}(x)C_{k0}(x)-
\tilde{R}_{ni,k1}(x)C_{k1}(x)\big),                          \label{1.6.19}
\end{equation}
\begin{equation}
R_{ni,lj}(x)-\tilde{R}_{ni,lj} (x)+\displaystyle\sum_{k=1}^{\infty}
\big(\tilde{R}_{ni,k0}(x)R_{k0,lj}(x)
-\tilde{R}_{ni,k1}(x)R_{k1,lj}(x)\big) = 0,                   \label{1.6.20}
\end{equation}
where $i, j=0,1$ and $n, l \in \mathbb{N}.$
The series in \eqref{1.6.19} and \eqref{1.6.20} converge absolutely and uniformly with respect to $x \in [a_{1}, b_{1}].$ 

For each fixed $x\in [a_1,b_1],$  relation
\eqref{1.6.19} can be treated as a system of linear equations with respect
to $C_{ni}(x),\; n\in {\mathbb N},\; i=0,1.$ 
But the series therein converges only "with brackets", i.e. the terms in them cannot be dissociated.
For this reason, it is inconvenient to use
\eqref{1.6.19} for solving  Inverse problem~2. 
Further we transform \eqref{1.6.19} into a linear equation in the corresponding Banach space
of sequences, see formula \eqref{1.6.33} below.

Let $V$ be the set of indices $u=(n,i),\; n \ge 1,\; i=0,1.$ 
For each fixed $x \in [a_{1}, b_{1}],$ we define the vector
$$
\psi(x)=[\psi_u(x)]_{u \in V}= \left[\begin{array}{l}\psi_{n0}(x)\\
\psi_{n1}(x) \end{array}\right]_{n \in \mathbb N}
$$
by the formulae
$$
\left[\begin{array}{l}\psi_{n0}(x)\\ \psi_{n1}(x) \end{array}\right] :=
\left[ \begin{array}{ll}\chi_{n} & -\chi_{n}\\ 0 & 1 \end{array}\right]
\left[\begin{array}{l}C_{n0}(x)\\ C_{n1}(x) \end{array}\right],
\quad 
\chi_n :=\left\{ \begin{array}{ll} \xi_n^{-1},\quad & \xi_n\ne 0, \\
0,\quad & \xi_n=0, \end{array}\right. \quad n \in \mathbb N.
$$
We also consider the block matrix
$$
H (x) = [H_{u,v}(x)]_{u, v \in V} \;=
\left[ \begin{array}{ll}H_{n0,k0}(x) & H_{n0,k1}(x)\\
H_{n1,k0}(x) & H_{n1,k1}(x) \end{array}\right]_{n,k \in \mathbb N},
\quad u=(n,i),\;v=(k,j),
$$
determined in the following way:
$$
\left[ \begin{array}{ll}H_{n0,k0}(x) & H_{n0,k1}(x)\\
H_{n1,k0}(x) & H_{n1,k1}(x) \end{array}\right] :=
\left[ \begin{array}{ll}\chi_{n} & -\chi_{n}\\
0 & 1 \end{array}\right]
\left[ \begin{array}{ll}R_{n0,k0}(x) & R_{n0,k1}(x)\\
R_{n1,k0}(x) & R_{n1,k1}(x) \end{array}\right]
\left[ \begin{array}{ll}\xi_{k} & 1\\ 0 & -1 \end{array}\right].
$$
Analogously we define $\tilde{\psi}(x)$ and $\tilde{H}(x)$ by replacing  $C_{ni}(x)$ by $\tilde C_{ni}(x)$
and $R_{ni,kj}(x)$ by $\tilde R_{ni,kj}(x)$  in the previous definitions.

Let us consider the Banach space $\mathcal B$ of bounded sequences
$\alpha =[\alpha_u]_{u \in V}$ with the norm
$\|\alpha \|_{\mathcal B} := \displaystyle\sup_{u \in V}|\alpha_u|.$ 
It follows from (\ref{1.6.6}) that $\psi(x),\tilde\psi(x) \in \mathcal B.$ Using \eqref{1.6.7}, for $n, k \in {\mathbb N}$ and $i,j = 0,1$ we obtain the estimates
\begin{equation}
|H_{ni,kj}(x)| \le \frac{C\xi_k}{|\sqrt{\theta_{n1}}-\sqrt{\theta_{k1}}|+1}, \quad 
|\tilde{H}_{ni,kj}(x)| \le \frac{C\xi_k}{|\sqrt{\theta_{n1}}-\sqrt{\theta_{k1}}|+1}.    \label{1.6.30}
\end{equation}
Consider the linear operators on $\mathcal B$ associated with the matrices $H(x)$ and $\tilde H(x):$
$$H(x) \alpha := [ y_u]_{u \in V}, \; y_u := \sum_{v \in V} H_{u,v}(x) \alpha_{v},  \quad \tilde H(x) \alpha := [ \tilde y_u]_{u \in V}, \; \tilde y_u := \sum_{v \in V} \tilde H_{u,v}(x) \alpha_{v}.$$
Due to \eqref{1.6.30}, \eqref{ser2}, and $\{\xi_n\}_{n\in{\mathbb N}} \in l_2,$ for each fixed $x$ the operators $H(x)$ and $\tilde H(x)$ are linear bounded operators in the space $\mathcal B.$

It is easy to see from \eqref{1.6.19} that for each fixed
$x\in  [a_{1}, b_{1}]$ the vector $\psi(x)$ satisfies the equation
\begin{equation}
\tilde{\psi}(x)=(I+\tilde{H}(x)) \psi(x)                  \label{1.6.33}
\end{equation}
in the Banach space $\mathcal B,$ where $I$ is the identity operator. 
From \eqref{1.6.20} it follows that $I- H(x)$ is the inverse operator to  $I+\tilde H(x).$ The existence of the inverse operator means that equation (\ref{1.6.33}) is uniquely solvable.

Thus, we obtain the following algorithm for solving Inverse problem~2 in the case $a_1<b_1.$
\medskip

{\bf Algorithm 1.}
Let the function $M(\lambda)$ be given. \\
1) Find the sequences $\{\alpha_n\}_{n=1}^\infty$ and $\{ \lambda_{n1}\}_{n=1}^\infty$  as the residues and the poles of $M(\lambda),$ respectively. \\
2)  Choose any model boundary value problem $\tilde L_1(T).$ Construct $\tilde\psi(x)$ and $\tilde H(x)$ for $x \in [a_{1}, b_{1}].$ \\
3) Find $\psi(x)$ by solving equation (\ref{1.6.33}), $x \in  [a_{1}, b_{1}].$ \\
4) Find $C(x, \lambda),$  $x \in  [a_{1}, b_{1}],$ from \eqref{1.6.10}. \\
5) Calculate $q(x) = (C''(x, \lambda) + \lambda C(x, \lambda))C^{-1}(x, \lambda),$ $x \in  [a_{1}, b_{1}].$
\medskip

So, we can find the potential on $[a_1, b_1]$ by the method of spectral mappings  if $a_1 < b_1.$ The case $a_1=b_1$ will be treated in the next section. 
\medskip

{\bf Remark 2.} If \eqref{z restriction} does not hold, Algorithm 1 can be obtained as well.
In this case, we have no asymptotic formulae \eqref{wn precise} for individual weight numbers. Instead of this, formulae can be found for the sums of the weight numbers in the groups of asymptotically close eigenvalues; see \cite{specmap-group, star, matrix}, where the analogous situations occur. Further, one should apply the form of the method of spectral mappings developed in \cite{specmap-group}.

\section{Solution of Inverse problem~2: the case $a_1=b_1.$ Solution of Inverse problem~1}
\label{point}
Further, we use the recursive structure of $T:$ 
$$T=[a_1, b_1] \cup T_2$$ 
with the time scale $T_2$ consisting only of isolated points and segments, see \eqref{T_m} for the definition of $T_m.$ Assume that $N+M-\mu_1 > 1,$ then the time scale $T_2$ contains at least two points.
Let us provide several relations between some objects for $T$ and the analogous ones for $T_2.$

 Consider the solution $S_{2}(x, \lambda),$ $x \in T_2,$ of the Sturm--Liouville equation \eqref{1} on $T_2^{0^2}$ satisfying the initial conditions
$$S_{2}(a_{2}, \lambda)  = 0, \quad  S^\Delta_{2}(a_{2}, \lambda) =  1.$$
If $T_2^{0^2}$ is the empty set, then the time scale $T_2$ is the union of two isolated points $a_2$ and $a_3.$ In this case all values of $S_{2}(x, \lambda)$ are completely determined by the initial conditions. 

Denote ${\cal D}_0^2(\lambda) := S_2(b_{N+M}, \lambda).$  This function is the characteristic functions of the boundary value problem $L_0(T_2).$  

 We introduce the function $\Phi_2(x,\lambda),$ $x \in T_2,$ which is the solution of equation (\ref{1}) under the boundary conditions
\begin{equation*}\label{Phi m}
\Phi_2^\Delta(a_2,\lambda)=1,\quad \Phi_2(b_{N+M},\lambda)=0.\end{equation*}
We also consider $M_2(\lambda) := \Phi_2(a_2, \lambda)$ which is the Weyl function for the Sturm--Liouville boundary value problems $L_j(T_2),$ $j=0,1.$ 
It is easy to see that
\begin{equation} \label{rec M}
M_{2}(\lambda) = \frac{\Phi(a_{2}, \lambda)}{\Phi^\Delta(a_{2}, \lambda)}.
\end{equation} 
Clearly, the functions $\Phi_2(x,\lambda)$ and $M_{2}(\lambda)$ are analogues of the functions $\Phi(x,\lambda)$ and $M(\lambda),$ respectively.

Now we are in position to solve the local inverse problem in the case $a_1 = b_1,$  which consists in recovering the value  $q(a_1)$  given $M(\lambda).$ 
For this purpose we use the following relation obtained in \cite{asymp}: 
\begin{equation}{\cal D}_0^{2}(\lambda) = \alpha_{11}^1(\lambda) \Theta_0(\lambda) - \alpha^1_{12}(\lambda) \Theta_1(\lambda) = \Theta_0(\lambda) - (a_{2} - b_1)\Theta_1(\lambda).\label{comp D_0}\end{equation}
The conditions $a_1 = b_1$ and $a_1 \in T^{0^2}$ guarantee that  $N+M-\mu_1 > 1$ and all values participating in the formula exist. We also use the following asymptotic formula: 
\begin{equation}\frac{\Theta_0(\lambda)}{{\cal D}_0^{2}(\lambda)} = \left\{
\begin{array}{cc}
(a_{2} - a_1)^2 (q(a_1)-\lambda) - (a_{2} - a_1)\rho i + 1 + o(1), & a_{2} < b_{2}, \\[4mm] 
\displaystyle (a_{2} - a_1)^2 (q(a_1)-\lambda) + \frac{a_{2} - a_{1}}{a_{3} - a_{2}} + 1 + o(1), & a_2 = b_2,
\end{array} \quad 
\lambda \to -\infty. \right. \label{comp_q}\end{equation} 
Note that formula \eqref{comp_q} can not be proved under the weaker assumptions on the potential than $q \in W^1_2[a_{l_k}, b_{l_k}],$ $k=\overline{1, N}.$ 
This is due to the fact that information about the values of $q$ in the isolated points can not be extracted from the leading terms in the polynomials $\beta^{l}_{ij}(\lambda).$ The details of the proof can be found in \cite{asymp}.

Taking~\eqref{comp_q} into account, we get the following algorithm for solving Inverse problem~2  in the case $a_1=b_1.$
\medskip

{\bf Algorithm 2.}
Let the function $M(\lambda)$ be given. \\
1) Construct $\Theta_0(\lambda)$ and $\Theta_1(\lambda).$ Find ${\cal D}_0^{2}(\lambda)$ with formula \eqref{comp D_0}. \\
2) Find $q(a_m)$ from \eqref{comp_q}.
\medskip

Algorithms~1 and~2 give the complete solution of the local inverse problem. Now we are ready to formulate the recursive algorithm for solving Inverse problem~1. Assume that conditions  \eqref{commensurability} and 
\eqref{condition on z} are fulfilled. 
\medskip

{\bf Algorithm 3.} Given the Weyl function $M(\lambda).$ \\
1) Construct $q(x)$ on $[a_{1}, b_{1}]$ using  Algorithm~1  or Algorithm~2. \\
2) If $T^{0^2} = [a_1, b_1],$ then terminate the algorithm. \\
3) Calculate $C^{\Delta^\nu}(b_1,\lambda)$ and $S^{\Delta^\nu}(b_1,\lambda)$ for $\nu=0,1.$
\\
4) Find $\Phi^{\Delta^\nu}(b_1,\lambda)$ for  $\nu=0,1$ by \eqref{Phi}.\\
5) Compute $\Phi(a_{2},\lambda)$ and $\Phi^\Delta(a_{2},\lambda)$ via jump conditions \eqref{jump conditions}.\\
6) Calculate  $M_{2}(\lambda)$ via \eqref{rec M}. Apply Algorithm~3 to $T_2$ given the Weyl function $M_2(\lambda).$
\medskip  

On step 6), we  run the Algorithm~3 to the time scale $T_2$ given the corresponding data $M_2(\lambda).$ 
Actually, we should re-designate $T=T_2$ and consider all other objects for this new $T$($=T_2$).
Then Algorithm~3 is repeated with the same notations. Call every its launch an iteration: the initial launch is the first iteration and so on. 

Conditions \eqref{commensurability} and \eqref{condition on z} formulated in the initial terms guarantee that the asymptotic formulae for the eigenvalues and the weight numbers  necessary for Algorithm~1 hold on every iteration as soon as the corresponding new $T$ begins with a segment. However, condition \eqref{condition on z} can be eliminated by the way mentioned in Remark 2.

After the $k$-th iteration, the potential is found on the set $\bigcup_{l=1}^k [a_l, b_l].$ Thus, the reconstruction of the potential will be fully completed after a finite number of iterations.\\


\end{document}